\documentclass[12pt]{amsart}

\usepackage{epsf,amsmath,amssymb,latexsym,amsthm}

\newcommand{\fix}{{\rm fix}}
\newcommand{\csum}{{\rm csum}}
\newcommand{\exc}{{\rm exc}}

\newcommand{\sumlim}{\sum\limits}
\newcommand{\prodlim}{\prod\limits}

\newtheorem{thm}{Theorem}[section]

\newtheorem{lem}[thm]{Lemma}

\newtheorem{cor}[thm]{Corollary}

\newtheorem{obs}[thm]{Observation}

\newtheorem{exa}[thm]{Example}

\theoremstyle{definition}
\newtheorem{defn}[thm]{Definition}

\newtheorem{prop}[thm]{Proposition}
\newtheorem{conj}[thm]{Conjecture}
\newtheorem{clm}[thm]{Claim}
\newcommand{\een}{\end{enumerate}}
\newcommand{\blem}{\begin{lem}}
\newcommand{\elem}{\end{lem}}
\newcommand{\bcl}{\begin{clm}}
\newcommand{\ecl}{\end{clm}}
\newcommand{\bthm}{\begin{thm}}
\newcommand{\ethm}{\end{thm}}
\newcommand{\bpr}{\begin{prop}}
\newcommand{\epr}{\end{prop}}
\newcommand{\bco}{\begin{cor}}
\newcommand{\eco}{\end{cor}}
\newcommand{\bcon}{\begin{conj}}
\newcommand{\econ}{\end{conj}}
\newcommand{\bde}{\begin{defn}}
\newcommand{\ede}{\end{defn}}
\newcommand{\bex}{\begin{exa}}
\newcommand{\eex}{\end{exa}}
\newcommand{\bobs}{\begin{obs}}
\newcommand{\eobs}{\end{obs}}
\newcommand{\bexe}{\begin{exe}}
\newcommand{\eexe}{\end{exe}}

\newcommand{\Z}{{\Bbb Z}}

\newcommand{\grn}{G_{r,n}}
\newcommand{\mulgrn}{G_{r_1,\dots,r_k;n}}

\begin{document}

\title[Statistics on the multi-colored permutation groups]{Statistics on the multi-colored permutation groups \\ }

\author{Eli Bagno}
\address{The Jerusalem College of Technology, Jerusalem, Israel}
\email{bagnoe@jct.ac.il}

\author{Ayelet Butman}
\address {Department of Computer Science, Faculty of Sciences, Holon Institute of Technology, PO Box
305, 58102 Holon, Israel} \email{ayeletb@hit.ac.il}

\author{David Garber}
\address{Department of Applied Mathematics, Faculty of Sciences, Holon Institute of Technology, PO Box 305,
58102 Holon, Israel} \email{garber@hit.ac.il}

\date{\today}

\maketitle
\begin{abstract}
We define an excedance number for the {\it multi-colored permutation
group} i.e. the wreath product $\mathbb{Z}_{r_1} \times \cdots
\times \mathbb{Z}_{r_k} \wr S_n$ and calculate its
multi-distribution with some natural parameters.

We also compute the multi-distribution of the parameters $\exc(\pi)$
and $\fix(\pi)$ over the sets of involutions in the multi-colored
permutation group. Using this, we count the number of involutions in
this group having a fixed number of excedances and absolute fixed
points.
\end{abstract}

\bibliographystyle{is-alpha}

\section{Introduction}

Let $r_1,\dots,r_k$ and $n$ be positive integers. The {\it
multi-colored permutation group} $G_{r_1,\dots,r_k;n}$ is the wreath
product:
$$(\mathbb{Z}_{r_1} \times \mathbb{Z}_{r_2} \times \cdots \times
\mathbb{Z}_{r_k}) \wr S_n.$$

The symmetric group $S_n$ is a special case for $r_i=1, 1 \leq i
\leq k$. In $S_n$ one can define the following well-known
parameters: Given $\sigma \in S_n$, $i \in [n]$ is {\it an excedance
of $\sigma$} if $\sigma(i)>i$. The number of excedances is denoted
by ${\rm exc}(\sigma)$. Two other natural parameters on $S_n$ are
the number of fixed points and the number of cycles of $\sigma$,
denoted by ${\rm fix}(\sigma)$ and ${\rm cyc}(\sigma)$ respectively.

Consider the following generating function over $S_n$:

$$P_n(q,t,s)=\sumlim_{\sigma \in
S_n}{q^{{\rm exc}(\sigma)}t^{{\rm fix}(\sigma)}s^{{\rm cyc}(\sigma)}}.$$

$P_n(q,1,1)$ is the classical Eulerian polynomial, while
$P_n(q,0,1)$ is the counter part for the derangements, i.e. the
permutations without fixed points, see \cite{Sta}.

In the case $s=-1$, the two polynomials $P_n(q,1,-1)$ and $P_n(q,0,-1)$
have simple closed formulas:

\begin{equation} \label{q=1}
P_n(q,1,-1)=-(q-1)^{n-1},
\end{equation}

\begin{equation}\label{q=-1}
P_n(q,0,-1)=-q[n-1]_q.
\end{equation}

Recently, Ksavrelof and Zeng \cite{KZ} proved some new recursive
formulas which induce the above equations. In \cite{BG}, the
corresponding excedance number for the colored permutation groups
$G_{r,n}=\mathbb{Z}_{r} \wr S_n$ was defined. It was proved there
that:
$$P_{\grn}(q,1,-1)=(q^r-1)P_{G_{r,n-1}}(q,1,-1),$$
$$P_{\grn}(q,0,-1)=[r]_q(P_{G_{r,n-1}}(q,0,-1) -q^{n-1}
[r]_q^{n-1}),$$ where $[n]_q=\frac{q^n-1}{q-1}$ and hence,
$$P_{\grn}(q,1,-1)=-\frac{(q^r-1)^n}{q-1},$$
$$P_{\grn}(q,0,-1)=-q [r]_q^n [n-1]_q.$$

\medskip
In this paper we generalize our parameters and formulas to the case
of the multi-colored permutation groups. Explicitly, denote $r=r_1
\cdots r_k$. We get the following theorems:

\begin{thm}\label{main}
$$P_{\mulgrn}(q,1,-1) =(q^r -1) P_{G_{r_1,\dots
,r_k;n-1}}(q,1,-1).$$ Hence,
$$P_{G_{r_1,\dots,r_k;n}}(q,1,-1)=\left( -1-K(q)\right)(q^r-1)^{n-1},$$
where  $$K(q)=K(q;r_1,\dots,r_k)=\sumlim_{m=1}^k r_{m+1} \cdots r_k
\sumlim_{t=1}^{r_m-1} q^{t \frac{r}{r_m}}.$$
\end{thm}

For the derangements, we have:

\begin{thm}\label{derange}
$$P_{\mulgrn} (q,0,-1) =\left(1+K(q) \right)
\left(P_{G_{r_1,\dots,r_k;n-1}}(q,0,-1) - \left( q^r+ K(q)
\right)^{n-1} \right).$$ Hence, we have:
\begin{tiny}
$$P_{G_{r_1,\dots,r_k;n}}(q,0,-1) = \left( q^r+ K(q) \right)
\left(1+K(q) \right)\cdot \left(\left(1+K(q) \right)^{n-2} -
\sumlim_{k=1}^{n-2} \left( q^r+ K(q) \right)^k \left(1+K(q)
\right)^{n-2-k}\right)$$
\end{tiny}
for all $n \geq 2$.
\end{thm}

An element $\sigma$ in $\mulgrn$  is called an {\it involution} if
$\sigma ^2=1$. The set of involutions in $\mulgrn$ will be denoted
by $I_{r_1,\dots,r_k;n}$.

In \cite{BGM}, the multi-distribution of the parameters ${\rm exc}$,
${\rm fix}$ and ${\rm csum}$ on the set of involutions in the
complex reflection groups was considered.  We cite the following
result from there. (The relevant definitions will be given in
Section \ref{involutions}).

\begin{thm} {\rm (See Corollary 5.2 in \cite{BGM}) }\\
The polynomial $\sumlim_{\pi \in \grn}u^{{\rm fix}(\pi)}v^{{\rm
exc_A}(\pi)}w^{{\rm \csum}(\pi)}$ is given by

\begin{equation}\label{pol}\sum_{j=n/2}^n (n-j)!\binom{n}{n-j,n-j,2j-n}\frac{u^{2j-n}(v+(r-1)w^r)^{n-j}}{2^{n-j}}\mu_r^{2j-n}.\end{equation}
where $\mu_r=1$ if $r$ is odd, and $\mu_r=1+w^{\frac{r}{2}}$ otherwise.
\end{thm}

Here, we generalize this result to $\mulgrn$. We prove:
\begin{thm}\label{invo}
The polynomial $\sumlim_{\pi \in \mulgrn}u^{{\rm fix}(\pi)}v^{{\rm
exc_A}(\pi)}w^{{\rm \csum}(\pi)}$ is given by

\begin{equation}\sum_{j=n/2}^n (n-j)!\binom{n}{n-j,n-j,2j-n}\frac{u^{2j-n}(v+(r-1)w^r)^{n-j}}{2^{n-j}}\mu^{2j-n}.\end{equation}
where $\mu=1$ if $r$ is odd, and $\mu=1+2^{\epsilon}
w^{\frac{r}{2}}$ otherwise (where\break $\epsilon = \# \{ r_i\ |\ 1
\leq i \leq k, r_i \equiv 0 \pmod 2 \}$). Hence, we have that the
number of involutions  $\pi \in \mulgrn$ with $\exc(\pi)=m$ is:
\begin{equation*}\left\{\begin{array}{cc}
y!\binom{n}{y,\ y,\ n-2y}{(\frac{r}{2})}^y & \qquad r \equiv 1 \pmod 2 \\
\sumlim_{j=\frac{n}{2}}^{n}{(n-j)!\binom{n}{n-j,\ n-j,\ j-y,\
y-n+j}{(\frac{r}{2})}^{n-j}2^{\epsilon(y-n+j)}} & \qquad r \equiv 0
\pmod 2
\end{array} \right.
\end{equation*}
where $y=\frac{m}{r}$.
\end{thm}

Note that every Abelian group $G$ can be presented as a direct
product of cyclic groups, and thus this work generalizes the
well-known excedance number to the wreath product of $S_n$ by any
Abelian group. Nevertheless, this parameter depends on the order of
the cyclic factors chosen to appear in the presentation of $G$.
Hence, it is an invariant of the pair $(G,(r_1,\dots,r_k))$ where $G
= \mathbb{Z}_{r_1} \times \cdots \times \mathbb{Z}_{r_k}$.

\medskip

The paper is organized as follows. In Section \ref{pre}, we give the
needed definitions. In Section \ref{stat}, we define the statistics
on $\mulgrn$. Section \ref{proof theorem 1} deals with the proof of
Theorem \ref{main}. Section \ref{derangement} deals with
derangements in $\mulgrn$ and the proof of Theorem \ref{derange}. In
Section \ref{involutions}, we deal with the set of involutions in
$\mulgrn$ and the proof of Theorem \ref{invo}.

\section{The group of multi-colored
permutations}\label{grn}\label{pre}

\bde Let $r_1,\dots,r_k$ and $n$ be positive integers. {\it The
group of multi-colored permutations of $n$ digits} is the wreath
product
$$G_{r_1,\dots,r_k;n}=(\mathbb{Z}_{r_1} \times \mathbb{Z}_{r_2} \times \cdots \times \mathbb{Z}_{r_k})\wr S_n
=(\mathbb{Z}_{r_1} \times \mathbb{Z}_{r_2} \times \cdots \times
\mathbb{Z}_{r_k})^n \rtimes S_n,$$ consisting of all the pairs
$(Z,\tau)$ where $Z=(z_i^j)$ ($1 \leq i \leq n,1 \leq j \leq k$) is
an $n \times k$ matrix such that the elements of column $j$ ($1 \leq
j \leq k$) belong to $\mathbb{Z}_{r_j}$ and $\tau \in S_n$. The
multiplication is defined by the following rule: Let $Z,U$ be two $n
\times k$ matrices as above and let $\sigma, \tau \in S_n$. Then

$$(Z,\tau) \cdot (U,\sigma)=((z_{i}^j+
u_{\tau^{-1}(i)}^j),\tau \circ \tau')$$ (here, in each column $j$,
the $+$ is taken modulo $r_j$). \ede

\bex Let $r_1=3,r_2=2,r_3=2,r_4=3$ and $n=3$.  Define
$$\pi_1=(Z_1,\tau_1)=\left(\ \left(\begin{array}{cccc}
0 & 1 & 0 & 2 \\
2 & 0 & 1 & 2 \\
1 & 1 & 0 & 1
\end{array}\right),\left(\begin{array}{ccc}
1 & 2 & 3\\
3 & 2 & 1
\end{array}\right)\ \right)$$
and
$$\pi_2=(Z_2,\tau_2)=\left(\ \left(\begin{array}{cccc}
0 & 0 & 1 & 0 \\
0 & 1 & 1 & 1 \\
2 & 1 & 0 & 2
\end{array}\right),\left(\begin{array}{ccc}
1 & 2 & 3\\
2 & 3 & 1
\end{array}\right)\ \right).$$

Then we have:
$$\pi_1 \cdot \pi_2 = \left(\ \left(\begin{array}{cccc}
2 & 0 & 0 & 2 \\
1 & 0 & 1 & 2 \\
2 & 0 & 0 & 1
\end{array}\right),\left(\begin{array}{ccc}
1 & 2 & 3\\
2 & 1 & 3
\end{array}\right)\ \right)$$

$$\pi_2 \cdot \pi_1=\left(\ \left(\begin{array}{cccc}
2 & 0 & 0 & 1 \\
2 & 1 & 0 & 0 \\
1 & 1 & 1 & 1
\end{array}\right),\left(\begin{array}{ccc}
1 & 2 & 3\\
1 & 3 & 2
\end{array}\right)\ \right).$$
\eex

Here is another description of the group $\mulgrn$. Consider the
alphabet
$$\Sigma=\{i^{[z_i^1,\dots,z_i^k]}\mid z_i^j \in \mathbb{Z}_{r_j},1
\leq i \leq n,1\leq j \leq k\}.$$ The set $\Sigma$ can be seen as
the set $[n]=\{1,\dots,n\}$, colored by $k$ palettes of colors, the
palette numbered $j$ having $r_j$ colors.

 If we denote by $\theta_j$ the cyclic operator which
colors the digit $i$ by first color from the $j$-th palette, then an
element of $\mulgrn$ is a {\it multi-colored permutation}, i.e. a
bijection $\pi: \Sigma \rightarrow \Sigma$ such that
$$\pi((\theta_1^{\epsilon_1} \circ \theta_2^{\epsilon_2} \circ \cdots \circ \theta_k^{\epsilon_k})(i))=
(\theta_1^{\epsilon_1} \circ \theta_2^{\epsilon_2} \circ \cdots
\circ \theta_k^{\epsilon_k})(\pi(i))$$
  where $\epsilon_i \in \{0,1 \}, 1 \leq i \leq k$.

In particular, if $k=1$ we get the group $G_{r_1,n}=C_{r_1} \wr
S_n$. This case has several subcases, for example if we take
$r=r_1=1$, then we get the symmetric group $S_n$, while $r=2$ yields
the hyperoctahedral group $B_n$, i.e.,  the classical Coxeter group
of type $B$.

\medskip

Here is an algebraic description of $\mulgrn$. Define the following
set of generators: $T=\{t_1,t_2,\dots,t_k,s_1,\dots,s_{n-1}\}$ with
the following relations:
\begin{itemize}

\item $t_i^{r_i}=1, (i \in \{1,\dots, k\})$

\item $(t_i s_1)^{2r_i}=1,(i \in \{1,\dots, k\})$

\item $s_i^2=1,(i \in \{1,\dots, n-1\})$

\item $s_is_js_i=s_js_is_j, (1 \leq i <j <n, j-i=1)$

\item $s_is_j=s_js_i, (1 \leq i <j<n, j-i>1)$

\item $t_i s_j=s_jt_i, (1 \leq i \leq k, 1<j<n).$

\end{itemize}

\begin{figure}[!ht]
\epsfxsize=10cm \epsfbox{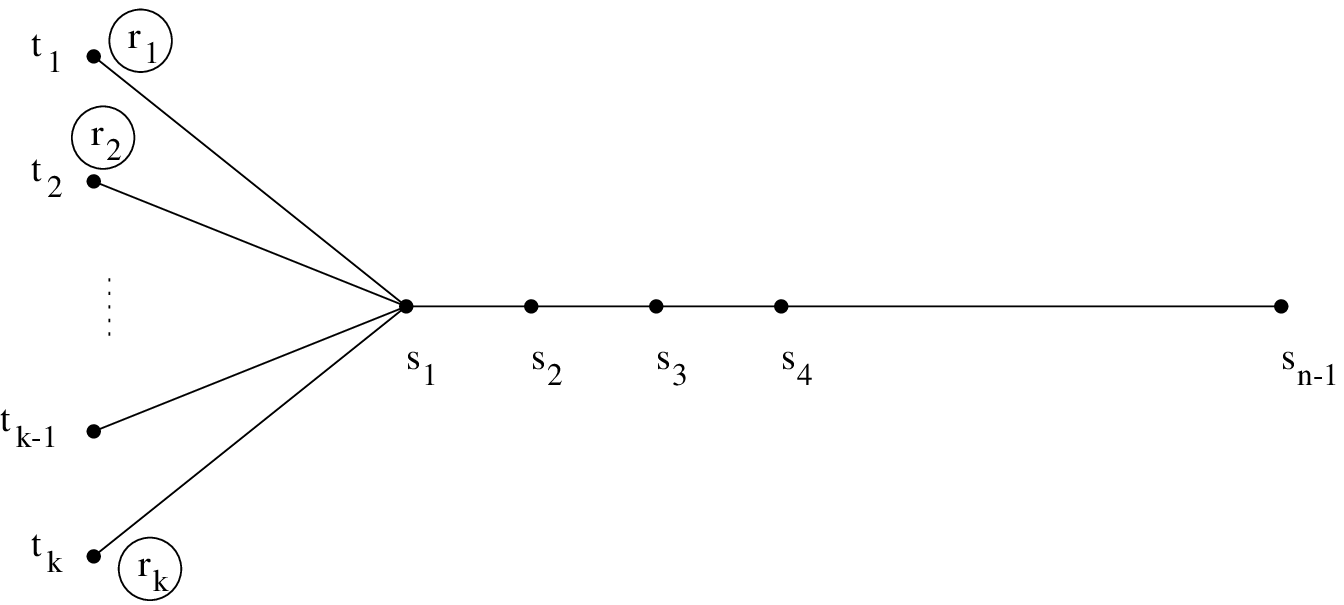} \caption{The "Dynkin
diagram" of $\mulgrn$}\label{muldyn}
\end{figure}

Realizing $t_i$ ($1 \leq i \leq k)$ as the multi-colored permutation
taking $1$ to $1^{\vec{e_i}}$ (where $\vec{e}_i$ is the $i$-th
standard vector) fixing pointwise the other digits, and $s_i$ as the
adjacent Coxeter transposition $(i,i+1)$ $(1 \leq i<n)$, it is easy
to see that $\mulgrn$ is actually the group generated by $T$ subject
to the above relations. A Dynkin-type diagram for $\mulgrn$ is
presented in Figure \ref{muldyn}.

\section{Statistics on $\mulgrn$}\label{stat}

We start by defining an order on the set:
$$\Sigma=\{i^{(z_i^1,\dots,z_i^k)}\mid z_i^j \in \mathbb{Z}_{r_k},1
\leq i \leq n,1\leq j \leq k\}.$$

Define $r_{\max}=\max \{r_1,\dots,r_k\}$. For any two vectors
$$\vec{v}=(v_1,\dots,v_k),\vec{w}=(w_1,\dots,w_k) \in \Z_{r_1}\times
\cdots \times \Z_{r_k},$$ we write $\vec{v}\prec \vec{w}$ if
$$
w_1 \cdot r_{\max}^{k-1}+\cdots +w_{k-1} \cdot  r_{\max}+w_k<v_1
\cdot r_{\max}^{k-1}+\cdots +v_{k-1} \cdot  r_{\max}+v_k.$$

For example, if $r_1=r_2=r_3=3$ then $(2,0,1) \prec (1,1,0)$.

 We also write $i^{\vec{v}} \prec j^{\vec{w}}$ if:

\begin{enumerate}

\item $\vec{v} \neq \vec{w}$ and $\vec{v} \prec \vec{w}$, or

\item $\vec{v} = \vec{w}$ and $i<j$.
\end{enumerate}

Based on this order, we define the {\it excedance set} of a
permutation $\pi$ on $\Sigma$ :
$${\rm Exc}(\pi)=\{i \in \Sigma \mid \pi(i) \succ i\},$$
and the {\it excedance number} is defined to be ${\rm exc}(\pi)=|{
\rm Exc}(\pi)|$.

\medskip

For simplifying the computations, we define the excedance number in
a different way. The set $\Sigma$ can be divided into layers,
according to the palettes. Explicitly, for each $\vec{v} \in
\mathbb{Z}_{r_1} \times \cdots \times \mathbb{Z}_{r_k}$, define the
layer $\Sigma^{\vec{v}}= \{ 1^{\vec{v}}, \dots ,n^{\vec{v}} \}$. We
call the layer $\Sigma^{\vec{0}}$ the {\it principal part} of
$\Sigma$. We will show that ${\rm exc}(\pi)$ can be computed using
parameters defined only on $\Sigma^{\vec{0}}$.

Let
$\pi=(\sigma,(z_1^1,\dots,z_1^k),(z_2^1,\dots,z_2^k),\dots,(z_n^1,\dots,z_n^k))
\in \mulgrn$ and let $1\leq p \leq k$. Define:
$${\rm csum}_p(\pi)=\sumlim_{i=1}^{n}z_i^p \cdot
\prodlim_{t=1}^{p-1}\chi (z_i^t=0),$$ where $\chi(P)$ is $1$ if the
property $P$ holds and $0$ otherwise.

The parameter ${\rm csum}_p(\pi)$ sums the colors of palette $p$
where a color of a digit is counted only if there are no colors of
preceding palettes on this digit.

Here is an easier way to understand the parameters ${\rm
csum}_p(\pi)$:

For
$\pi=(\sigma,(z_1^1,\dots,z_1^k),(z_2^1,\dots,z_2^k),\dots,(z_n^1,\dots,z_n^k))
\in \mulgrn$, write the $n \times k$ matrix $Z=(z_i^j)$. Then,
${\rm csum}_p$ is just the sum of the elements of the $p$-th
column where we are ignoring the elements which are not leading in
their rows.

\begin{exa} Let $$\pi=\left(\begin{pmatrix} 1 & 2 & 3  \\
                                      3 & 1 & 2
\end{pmatrix},(1,2,0,1),(0,0,1,2),(0,2,1,1)\right) \in G_{2,3,2,3;3}.$$

Then $Z=\begin{pmatrix} \mathbf{1} & 2 & 0 & 1\\
                        0 & 0 & \mathbf{1} & 2 \\
                        0 & \mathbf{2} & 1 & 1\end{pmatrix}$

and thus we have:

$${\rm csum}_1(\pi)=1,{\rm csum}_2(\pi)=2,{\rm csum}_3(\pi)=1,{\rm
csum}_4(\pi)=0.$$

\end{exa}

Now define:

\centerline{${\rm Exc}_A(\pi)=\{i \in [n-1] \mid \pi(i) \succ i\}$
and ${\rm exc}_A(\pi)=|{\rm Exc}_A(\pi)|$.}

\medskip

\bpr Let $\pi=(Z,\sigma)$. Write $r=\prodlim_{j=1}^{k}r_j$. Then:

$${\rm exc}(\pi)=r \cdot {\rm
exc}_A(\pi)+\sumlim_{p=1}^{k} \left( {\rm csum}_p(\pi)\cdot
\prodlim_{q=1,q\neq p}^{k}r_q \right).$$ \epr

\begin{proof}
Let $i \in [n]$. Write $\pi(i^{\vec{0}})=j^{\vec{z}_i}$. We divide
our treatment according to $\vec{z}_i =(z_i^1,\dots,z_i^k)$.

\begin{itemize}
\item $\vec{z}_i=\vec{0}$: In this case, $i \in {\rm Exc}_A (\pi)$
if and only if $\sigma(i)>i$ or in other words:
$\pi\left(i^{\vec{0}}\right) \succ i^{\vec{0}}$. This happens, if
and only if, for each $\vec{\alpha}=(\alpha^1,\dots,\alpha^k)$ where
$0 \leq \alpha^t \leq r_t-1$, we have
$\pi\left(i^{\vec{\alpha}}\right) \succ i^{\vec{\alpha}}$. Thus $i$
contributes $\prodlim_{j=1}^{k}r_j=r$ to ${\rm exc}(\pi)$.

\item $\vec{z}_i=(z_i^1,\dots,z_i^k) \neq \vec{0}$. In this case $i=i^{\vec{0}} \not \in {\rm
Exc}( \pi)$. We check now for which $\vec{v}$, $i^{\vec{v}} \in \rm
Exc(\pi)$. Since $\pi(i^{\vec{0}})=j^{\vec{z_i}}$, we have
$\pi(i^{\vec{v}})=j^{\vec{v}+\vec{z_i}}$.
 Let $m \in
\{1,\dots, k\}$ be the minimal index such that $z_i^m \neq 0$ and
$z_i^t=0$ for all $t<m$. Note that for all $\vec{0} \succ \vec{v}
\succ (0,\dots,0,r_m-z_i^m,0,\dots ,0)$,
$\pi(i^{\vec{v}})=j^{\vec{v}+\vec{z_i}} \prec i^{\vec{v}}$, hence
$i^{\vec{v}} \not \in {\rm Exc}(\pi)$.

Now, for all $$\hspace{20pt}(0,\dots,0,r_m-z_i^m,0,\dots ,0) \succ
\vec{v} \succ (0,\dots,0, r_m-1,r_{m+1}-1,\dots,r_k-1),$$
$\pi(i^{\vec{v}})=j^{\vec{v}+\vec{z_i}} \succ i^{\vec{v}}$, and
hence $i^{\vec{v}} \in {\rm Exc}(\pi)$. So, it contributes $z_i^m
\cdot r_{m+1} \cdots r_k$ elements to the excedance set.

\medskip

In the same way, for each $\vec{w}=(\alpha_1,\dots,
\alpha_{m-1},0,\dots,0)\neq \vec{0}$ and for all
$$\hspace{20pt}(0,\dots,0,r_m-z_i^m,0,\dots ,0) \succ \vec{v} \succ (0,\dots,0,
r_m-1,r_{m+1}-1,\dots,r_k-1),$$
$\pi(i^{\vec{w}+\vec{v}})=j^{\vec{w}+\vec{v}+\vec{z_i}} \succ
i^{\vec{w}+\vec{v}}$, and hence $i^{\vec{w}+\vec{v}} \in {\rm
Exc}(\pi)$. So it contributes $(r_1 \cdots r_{m-1} -1) \cdot z_i^m
\cdot r_{m+1} \cdots r_k$ elements to the excedance set.

Hence, this $i$ contributes $$r_1 \cdots r_{m-1} \cdot z_i^m \cdot
r_{m+1} \cdots r_k = z_i^m \prodlim_{q=1,q\neq m}^{k} r_q.$$

\end{itemize}

\medskip

Now, we sum the contributions over all $i \in \{1,2,\dots,n\}$.
Since we have ${\rm exc}_A(\pi)$ digits which satisfy
$\vec{z}_i=\vec{0}$ and $\sigma(i)>i$, their total contribution is
$r \cdot {\rm exc}_A(\pi)$, which is the first summand of ${\rm
exc}(\pi)$.

The other digits have $\vec{z}_i \neq \vec{0}$, so their
contribution is
$$\sumlim_{\{i\mid \vec{z}_i \neq \vec{0}\}} \left( z_i^m \prodlim_{q=1,q\neq m}^{k}
r_q \right) =\sumlim_{p=1}^{k} \left( \left( \sumlim_{i=1}^{n}z_i^p
\cdot \prodlim_{t=1}^{p-1}\chi (z_i^t=0) \right) \prodlim_{q=1,q\neq
p}^{k}r_q \right)=$$
$$=\sumlim_{p=1}^{k} \left( {\rm
csum}_p(\pi)\cdot \prodlim_{q=1,q\neq p}^{k}r_q \right),$$ which
is the second summand of ${\rm exc}(\pi)$, and hence we are done.
\end{proof}

\begin{exa}
Let
$$\pi =\left( \begin{array}{ccc} 1 & 2 & 3 \\ 3^{(0,0)} & 1^{(2,1)} & 2^{(0,1)}  \end{array} \right) = \left(\ \left(\begin{array}{cc}
0 & 0  \\
\mathbf{2} & 1 \\
0 & \mathbf{1}
\end{array}\right),\left(\begin{array}{ccc}
1 & 2 & 3\\
3 & 1 & 2
\end{array}\right)\ \right) \in G_{3,2;3}.$$
We write $\pi$ in its extended form:

$$\hspace{50pt}\left. \begin{array}{|ccc|ccc|||ccc} \checkmark & \nabla & & \checkmark& &\maltese & \checkmark& & \\
1^{(1,0)} & 2^{(1,0)} & 3^{(1,0)} & 1^{(0,1)} & 2^{(0,1)} &
3^{(0,1)} & 1^{(0,0)} &
2^{(0,0)} & 3^{(0,0)} \\
3^{(1,0)} & 1^{(0,1)} & 2^{(1,1)} & 3^{(0,1)} & 1^{(2,0)} &
2^{(0,0)} & 3^{(0,0)} & 1^{(2,1)} & 2^{(0,1)}
\end{array} \right) $$

$$\hspace{-50pt}\left( \begin{array}{ccc|ccc|ccc|}\checkmark & \nabla & \maltese & \checkmark& \nabla & & \checkmark& \nabla & \maltese \\
 1^{(2,1)} & 2^{(2,1)} & 3^{(2,1)} & 1^{(2,0)} & 2^{(2,0)} & 3^{(2,0)} &
1^{(1,1)} & 2^{(1,1)} & 3^{(1,1)} \\
3^{(2,1)} & 1^{(1,0)} & 2^{(2,0)} & 3^{(2,0)} & 1^{(1,1)} &
2^{(2,1)} & 3^{(1,1)} & 1^{(0,0)} & 2^{(1,0)}
\end{array} \right. $$

We have ${\rm exc}(\pi)=13$, while ${\rm csum}_1 (\pi)=2$ and
${\rm csum}_2(\pi)=1$.

\end{exa}

\medskip

Recall that any permutation of $S_n$ can be decomposed into a
product of disjoint cycles. This notion can be easily generalized to
the group $\mulgrn$ as follows. Given any $\pi \in \mulgrn$ we
define the {\it cycle number} of $\pi=(Z,\sigma)$ to be the number
of cycles in $\sigma$.

We say that $i \in [n]$ is an {\it absolute fixed point} of $\pi \in
\mulgrn$ if $\sigma (i)=i$.

\section{Proof of Theorem \ref{main}}\label{proof theorem 1}

In this section we prove Theorem \ref{main}. The way to prove this
type of identities is to construct a subset $S$ of $\mulgrn$ whose
contribution to the generating function is exactly the right side of
the identity. Then, we have to construct a killing involution on
$\mulgrn - S$, i.e., an involution on $\mulgrn - S$ which preserves
the number of excedances but changes the sign of every element of
$\mulgrn -S$ and hence shows that $\mulgrn -S$ contributes nothing
to the generating function.

Recall that $r=r_1 \cdots r_k$. We divide $\mulgrn$ into $2r+1$
disjoint subsets as follows:

$$K=\{\pi \in \mulgrn \mid |\pi(n)|\neq n , |\pi(n-1)| \neq n\},$$

$$T^{\vec{v}}_{n}=\{\pi \in \mulgrn \mid \pi(n)=n^{\vec{v}}\}, \qquad
(\vec{v} \in \mathbb{Z}_{r_1} \times \cdots \times
\mathbb{Z}_{r_k}),$$

$$R^{\vec{v}}_{n}=\{\pi \in \mulgrn \mid \pi(n-1)=n^{\vec{v}}\}, \qquad
(\vec{v} \in \mathbb{Z}_{r_1} \times \cdots \times
\mathbb{Z}_{r_k}),$$

We first construct a killing involution on the set $K$. Let $\pi \in
K$. Define $\varphi: K \to K$ by
$$\pi ' = \varphi(\pi)=(\pi(n-1),\pi(n)) \pi.$$
Note that $\varphi$ exchanges $\pi(n-1)$ with $\pi(n)$. It is
obvious that $\varphi$ is indeed an involution.

We will show that $\rm{exc}(\pi)=\rm{exc}(\pi ')$. First, for
$i<n-1$, it is clear that $i \in {\rm Exc}(\pi)$ if and only if $i
\in {\rm Exc}(\pi')$. Now, as $\pi(n-1) \neq n$, $n-1 \notin {\rm
Exc}(\pi)$ and thus $n \notin {\rm Exc}(\pi')$. Finally, $\pi(n)\neq
n$ implies that $n-1 \notin {\rm Exc}(\pi')$ and thus ${\rm
exc}(\pi)={\rm exc}(\pi').$

On the other hand, $\rm{cyc}(\pi)$ and $\rm{cyc}(\pi')$ have
different parities due to a multiplication by a transposition.
Hence, $\varphi$ is indeed a killing involution on $K$.

\medskip

We turn now to the sets  $T^{\vec{v}}_{n} \quad
(\vec{v}=(z_n^1,\dots ,z_n^k) \in \mathbb{Z}_{r_1} \times \cdots
\times \mathbb{Z}_{r_k})$. Note that there is a natural bijection
between $T^{\vec{v}}_{n}$ and $G_{r_1,\dots,r_k;n-1}$ defined by
ignoring the last digit. Let $\pi \in T^{\vec{v}}_{n}$. Denote the
image of $\pi \in T^{\vec{v}}_{n}$ under this bijection by $\pi'$.
Since $n \not\in \rm{Exc}_A(\pi)$, we have
$\rm{exc}_A(\pi)=\rm{exc}_A(\pi')$.

Let $m \in \{1,\dots, k\}$ be the minimal index such that $z_n^m
\neq 0$ and $z_n^t=0$ for all $t<m$. Then, ${\rm
csum}_m(\pi')={\rm csum}_m(\pi)-z_n^m$, and ${\rm
csum}_p(\pi')=\rm{csum}_p(\pi)$ for $1 \leq p \leq k,p\neq m$.
Finally, since $n$ is an absolute fixed point of $\pi$,
$\rm{cyc}(\pi ')= \rm{cyc}(\pi)-1$. Hence, we get that the total
contribution of $T^{\vec{v}}_{n}$ is:
$$P_{T^{\vec{v}}_{n}} = -q^{z_n^m \prodlim_{q=1,q\neq m}^{k} r_q} P_{G_{r_1,\dots,r_k;n-1}}(q,1,-1) =
 -q^{z_n^m \frac{r}{r_m}} P_{G_{r_1,\dots,r_k;n-1}}(q,1,-1),$$
where $m$ is defined as above.

\medskip

Now, we treat the sets $R^{\vec{v}}_{n} \quad (\vec{v}=(z_n^1, \dots
,z_n^k) \in \mathbb{Z}_{r_1} \times \cdots \times
\mathbb{Z}_{r_k})$. There is a bijection between $R^{\vec{v}}_{n}$
and $T^{\vec{v}}_{n}$ using the same function $\varphi$ we used
above. Let $\pi \in R_{n}^{\vec{v}}$. Define $\varphi:
R^{\vec{v}}_{n} \to T^{\vec{v}}_{n}$ by
$$\pi ' = \varphi(\pi)=(\pi(n-1),\pi(n)) \pi.$$

When we compute the change in the excedance, we split our
treatment into two cases: $\vec{v}=\vec{0}$ and $\vec{v} \neq
\vec{0}$.


We start with the case $\vec{v}=\vec{0}$. Note that $n-1 \in
\rm{Exc}_A(\pi)$ (since $\pi(n-1)=n$) and $n \not\in
\rm{Exc}_A(\pi)$. On the other hand, in $\pi'$, $n-1,n \not\in
\rm{Exc}_A(\pi')$. Hence, $\rm{exc}_A (\pi)-1=\rm{exc}_A(\pi')$.

Now, for the case $\vec{v}\neq \vec{0}:$ $n-1,n \not\in
\rm{Exc}_A(\pi)$ (since $\pi(n-1)=n^{\vec{v}}$ is not an
excedance). We also have: $n-1,n \not\in \rm{Exc}_A(\pi')$ and
thus $\rm{Exc}_A(\pi)=\rm{Exc}_A(\pi ')$ for $\pi \in
R^{\vec{v}}_{n}$ where $\vec{v}\neq\vec{0}$.

In both cases, we have that ${\rm csum_p}(\pi) =\rm{csum}_p(\pi')$
for each $1\leq p\leq k$. Hence, we have that ${\rm
exc}(\pi)-r=\rm{exc}(\pi ')$ for $\vec{v}=\vec{0}$ and
$\rm{exc}(\pi)={\rm exc}(\pi ')$ for $\vec{v} \neq \vec{0}$.

As before, the number of cycles changes its parity due to the
multiplication by a transposition, and hence: $(-1)^{{\rm
cyc}(\pi)} = -(-1)^{\rm{cyc}(\pi')}$.

Hence, the total contribution of the elements in $R^{\vec{v}}_{n}$
is
$$q^r P_{G_{r_1,\dots,r_k;n-1}}(q,1,-1)$$ for $\vec{v}=\vec{0}$, and
$$ q^{z_n^m \prodlim_{q=1,q \neq m}^k r_q} P_{G_{r_1,\dots,r_k;n-1}}(q,1,-1)=
q^{z_n^m \frac{r}{r_m}} P_{G_{r_1,\dots,r_k;n-1}}(q,1,-1)$$ for
$\vec{v} \neq \vec{0}$.

\bigskip

In order to calculate $\sumlim_{\vec{v} \in \mathbb{Z}_{r_1}
\times \cdots \times \mathbb{Z}_{r_k}} P_{T_n^{\vec{v}}}$ and
$\sumlim_{\vec{v} \in \mathbb{Z}_{r_1} \times \cdots \times
\mathbb{Z}_{r_k}} P_{R_n^{\vec{v}}}$, we have to divide
$\mathbb{Z}_{r_1} \times \cdots \times \mathbb{Z}_{r_k}$ into sets
according to the minimal index $m$ such that $z_n^m \neq 0$ and
$z_n^t=0$ for all $t<m$.

For each $m \in \{ 1, \dots, k+1 \}$, denote:
$$W_m = \{ \vec{v} = (z_n^1, \dots , z_n^k) \in \mathbb{Z}_{r_1} \times
\cdots \times \mathbb{Z}_{r_k} | z_n^m \neq 0, z_n^t=0 ,\forall t<m
\}.$$ Note that $W_{k+1} =\{\vec{0}\}$.

It is easy to see that $\{ W_1, \dots, W_k, W_{k+1} \}$ is a
partition of\break $\mathbb{Z}_{r_1} \times \cdots \times
\mathbb{Z}_{r_k}$.

Hence \begin{tiny}$$\sumlim_{\vec{v} \in \mathbb{Z}_{r_1} \times
\cdots \times \mathbb{Z}_{r_k}} P_{T_n^{\vec{v}}} =
\sumlim_{m=1}^{k+1} \left( \sumlim_{\vec{v} \in W_m}
P_{T_n^{\vec{v}}} \right) =\left( -1+ \sumlim_{m=1}^k r_{m+1}
\cdots r_k \sumlim_{t=1}^{r_m-1} -q^{t \frac{r}{r_m}}
\right)P_{G_{r_1,\dots,r_k;n-1}}(q,1,-1) .$$
\end{tiny}

Similarly, we get:
\begin{tiny}
$$\sumlim_{\vec{v} \in \mathbb{Z}_{r_1} \times \cdots
\times \mathbb{Z}_{r_k}} P_{R_n^{\vec{v}}} = \left( q^r+
\sumlim_{m=1}^k r_{m+1} \cdots r_k \sumlim_{t=1}^{r_m-1} q^{t
\frac{r}{r_m}} \right)P_{G_{r_1,\dots,r_k;n-1}}(q,1,-1).$$
\end{tiny}

 Now, if we sum up all the parts, we get:
\begin{eqnarray*}
\hspace{-40pt}P_{\mulgrn}(q,1,-1) &=&\sumlim_{\vec{v} \in
\mathbb{Z}_{r_1} \times \cdots \times \mathbb{Z}_{r_k}}
P_{T_n^{\vec{v}}}+ \sumlim_{\vec{v} \in \mathbb{Z}_{r_1} \times
\cdots \times \mathbb{Z}_{r_k}} P_{R_n^{\vec{v}}}=
\\ \hspace{-40pt}&=&(q^r -1) P_{G_{r_1,\dots ,r_k;n-1}}(q,1,-1)
\end{eqnarray*}
as needed.

\medskip

From now on, we denote
$$K(q)=K(q;r_1,\dots,r_k)=\sumlim_{m=1}^k r_{m+1} \cdots r_k
\sumlim_{t=1}^{r_m-1} q^{t \frac{r}{r_m}}.$$ Note that:
$$1+K(q)=\sumlim_{\vec{v} \in \mathbb{Z}_{r_1} \times \cdots \times \mathbb{Z}_{r_k}}q^{{\rm exc}(1^{\vec{v}})}.$$

\medskip

Now, for $n=1$, $G_{r_1,\dots,r_k;1}$ is $\mathbb{Z}_{r_1} \times
\cdots \times \mathbb{Z}_{r_k}$, and thus
$$P_{G_{r_1,\dots,r_k;1}}(q,1,-1)=-1-K(q).$$ Hence,
we have
$$P_{G_{r_1,\dots,r_k;n}}(q,1,-1)=\left( -1-K(q)\right)(q^r-1)^{n-1},$$
and we have finished the proof of Theorem \ref{main}.

\section{Derangements in $\mulgrn$ and the proof of Theorem \ref{derange}}\label{derangement}

We start with the definition of a derangement.

\bde An element $\sigma \in \mulgrn$ is called a {\it derangement}
if it has no absolute fixed points, i.e. $|\pi(i)|\neq i$ for
every $i \in [n]$. Denote by $D_{r_1,\dots,r_k;n}$ the set of all
derangements in $\mulgrn$. \ede

We prove now Theorem \ref{derange}.

We divide $D_{r_1,\dots,r_k;n}$ into $r+2=r_1 r_2 \cdots r_k +2$
disjoint subsets in the following way:

$$A^{\vec{v}}_{r_1,\dots,r_k;n}=\{\pi \in D_{r_1,\dots,r_k;n} \mid \pi(2)=1^{\vec{v}},|\pi(1)| \neq 2\}, \vec{v} \in \mathbb{Z}_{r_1} \times \cdots
\times \mathbb{Z}_{r_k}.$$

$$B_{r_1,\dots,r_k;n}=\{\pi \in D_{r_1,\dots,r_k;n} \mid |\pi|=(123 \cdots n)\}.$$

$$\hat{D}_{r_1,\dots,r_k;n}=D_{r_1,\dots,r_k;n} - \left(\bigcup_{\vec{v} \in \mathbb{Z}_{r_1} \times \cdots \times \mathbb{Z}_{r_k}}{A^{\vec{v}}_{r_1,\dots,r_k;n}} \cup B_{r_1,\dots,r_k;n}\right).$$

We start by constructing a killing involution $\varphi$ on
$\hat{D}_{r_1,\dots,r_k;n}$. Given any $\pi \in
\hat{D}_{r_1,\dots,r_k;n}$, let $i$ be the first number such that
$|\pi(i)| \neq i+1$. Define
$$\pi'=\varphi(\pi)=(\pi(i),\pi(i+1)) \pi.$$

It is easy to see that $\varphi$ is a well-defined involution on
$\hat{D}_{r_1,\dots,r_k;n}$. We proceed to prove that ${\rm
exc}(\pi)={\rm exc}(\pi')$. Indeed, ${\rm csum}_p(\pi)={\rm
csum}_p(\pi')$ for all $1 \leq p \leq k$.

Let $i$ be the first number such that $|\pi(i)| \neq i+1$, so that
in the pass from $\pi$ to $\pi'$ we exchange $\pi(i)$ with
$\pi(i+1)$. For every $j \neq i,i+1$, clearly $j \in {\rm
Exc}_A(\pi)$ if and only if $j \in {\rm Exc}_A(\pi')$. Since $\pi
\in D_{r_1,\dots,r_k;n}$, $|\pi(i)| \neq i+1$ and $|\pi(j)|=j+1$ for
$j<i$, we have that $|\pi(i)|,|\pi(i+1)| \in \{1,i+2,\dots,n\}$.
Thus, exchanging $\pi(i)$ with $\pi(i+1)$ does not change ${\rm
Exc}_A(\pi)$.

Note also that the parity of ${\rm cyc}(\pi')$ is opposite to the
parity of ${\rm cyc}(\pi)$ due to the multiplication by a
transposition. Hence, we have proven that $\varphi$ is indeed a
killing involution.

Now, let us calculate the contribution of each set in our
decomposition to $P_{\mulgrn}(q,0,-1)$. As we have shown,
$\hat{D}_{r_1,\dots,r_k;n}$ contributes nothing.

Let $\vec{v}=(z_2^1,\dots,z_2^n)$. Define a bijection
$$\psi: A^{\vec{v}}_{r_1,\dots,r_k;n} \to D_{r_1,\dots,r_k;n-1}$$ by: $\psi(\pi)=\pi'$
where $\pi'(1)=(|\pi(1)|-1)^{z_1(\pi)}$ and for $j>1$,
$\pi'(j)=(|\pi(j+1)|-1)^{z_{j+1}(\pi)}$. For example, if
$\pi=(3^{(0,1,0)}1^{(0,0,0)}4^{(2,2,2)}2^{(0,0,1)})$, then $\pi
'=(2^{(0,1,0)}3^{(2,2,2)}1^{(0,0,1)})$. It is easy to see that
${\rm exc}_A(\pi)={\rm exc}_A{(\pi')}$. On the other hand, let $m
\in \{1,\dots, k\}$ be the minimal index such that $z_2^m \neq 0$
and $z_n^t=0$ for all $t<m$.  Then ${\rm csum}_m(\pi')={\rm
csum}_m(\pi)-z_2^m$, and ${\rm csum}_p(\pi')={\rm csum}_p(\pi)$
for $1 \leq p \leq k,p\neq m$.

We have also: ${\rm cyc}(\pi)={\rm cyc}(\pi')$ and thus the
contribution of $A^{\vec{v}}_{r,n}$ to $P_{\mulgrn}(q,0,-1)$ is
$$P_{A_{r_1,\dots,r_k;n}^{\vec{v}}} = q^{z_2^m \prodlim_{q=1,q\neq m}^{k} r_q} P_{G_{r_1,\dots,r_k;n-1}}(q,0,-1) =
 q^{z_2^m \frac{r}{r_m}} P_{G_{r_1,\dots,r_k;n-1}}(q,0,-1)$$
where $m$ is defined as above.

 Finally, we treat the set $B_{r_1,\dots,r_k;n}$. Here for every $\pi
\in B_{r_1,\dots,r_k;n}$ we have ${\rm cyc}(\pi)=1$.

We calculate now:
\begin{eqnarray*}
P_{B_{r_1,\dots,r_k;n}}(q,0,-1)&=&-\prodlim_{s=1}^{n}\sumlim_{\vec{v}_s
\in \mathbb{Z}_{r_1}\times \cdots \times
\mathbb{Z}_{r_k}}q^{{\rm{exc}}(2^{\vec{v}_1}3^{\vec{v}_2} \cdots
n^{\vec{v}_{n-1}} 1^{\vec{v}_n})}=\\
&=& - \left( q^r+ K(q) \right)^{n-1}\cdot \left(1+K(q) \right).
\end{eqnarray*}

 To summarize, we get:
\begin{eqnarray*}
P_{\mulgrn} (q,0,-1) & =&\left( \sumlim_{\vec{v}=(z_2^1,\dots,z_2^k)
\in \mathbb{Z}_{r_1} \times \cdots \times
\mathbb{Z}_{r_k}}{P_{A_{r_1,\dots,r_k;n}^{\vec{v}}}}\right)+P_{B_{r_1,\dots,r_k;n}}=\\
&=& \left( \sumlim_{m=1}^{k+1} \sumlim_{\vec{v} \in
W_m}{P_{A_{r_1,\dots,r_k;n}^{\vec{v}}}} \right)+P_{B_{r_1,\dots,r_k;n}}=\\
&=& \left(1+ K(q)\right)P_{G_{r_1,\dots,r_k;n-1}}(q,0,-1)-\\
& &-\left( q^r+ K(q) \right)^{n-1}\cdot
\left(1+K(q) \right)= \\
&=& \left(1+K(q) \right) \left(P_{G_{r_1,\dots,r_k;n-1}}(q,0,-1) -
\left( q^r+ K(q) \right)^{n-1} \right),
\end{eqnarray*}
so we get:
$$P_{\mulgrn} (q,0,-1) =\left(1+K(q) \right)
\left(P_{G_{r_1,\dots,r_k;n-1}}(q,0,-1) - \left( q^r+ K(q)
\right)^{n-1} \right).$$

Now, for $n=2$ we have:
$$P_{G_{r_1,\dots,r_k;2}}(q,0,-1) = \left( q^r+ K(q) \right)\cdot
\left(1+K(q) \right).$$

\medskip

By a direct computation, one can now get:
\begin{tiny}
$$P_{G_{r_1,\dots,r_k;n}}(q,0,-1) = \left( q^r+ K(q) \right)
\left(1+K(q) \right)\cdot \left(\left(1+K(q) \right)^{n-2} -
\sumlim_{k=1}^{n-2} \left( q^r+ K(q) \right)^k \left(1+K(q)
\right)^{n-2-k}\right)$$
\end{tiny}
for all $n \geq 2$ as in Theorem \ref{derange}.

\section{Involutions in $\mulgrn$}\label{involutions}

We recall that an element $\sigma$ in $\mulgrn$
is called an {\it involution} if
$\sigma ^2=1$. The set of involutions in $\mulgrn$ will be denoted
by $I_{r_1,\dots,r_k;n}$.

We consider the multi-distribution of the parameters ${\rm exc}$,
${\rm fix}$ and ${\rm csum}$ on $I_{r_1,\dots,r_k;n}$, where ${\rm
csum}(\pi)$ here is the total contribution of all the ${\rm
csum}_p$-s from all the palettes:
$${\rm csum}(\pi)= {\rm exc}(\pi) -r \cdot {\rm exc}_A(\pi) =\sumlim_{p=1}^{k} \left( {\rm csum}_p(\pi)\cdot
\prodlim_{q=1,q\neq p}^{k}r_q \right)   .$$

\medskip

We start by classifying the involutions of $\mulgrn$. As in the case
of $\grn$, each involution of $\mulgrn$ can be decomposed into a
product of 'atomic' involutions of two types: absolute fixed points
and $2$-cycles.

An absolute fixed point must be of the form $\pi(i)=i^{\vec{v}}$
where $2 \vec{v}=\vec{0}$.

The $2$-cycles have the form $\pi(i)= j^{\vec{v}_1}$; $\pi(j)=
i^{\vec{v}_2}$ where $\vec{v}_1,\vec{v}_2 \in \mathbb{Z}_1 \times
\cdots \times \mathbb{Z}_k$ and $\vec{v}_1+\vec{v}_2=\vec{0}$.



\medskip

Now, we compute recurrence and explicit formulas for
$$f_{r_1,\dots,r_k;n}(u,v,w) = \sum_{\pi\in
I_{r_1,\dots,r_k;n}}u^{\fix(\pi)}v^{\exc(\pi)}w^{\csum(\pi)}.$$

Let $\pi$ be any involution in $I_{r_1,\dots,r_k;n}$. Then we have
either $\pi(n)=n^{\vec{v}}$ or $\pi(n)=k^{\vec{v}}$ with $k<n$.

For $\pi \in I_{r_1,\dots, r_k;n}$  such that
$\pi(n)=n^{\vec{v}}$, define $\pi' \in I_{r_1,\dots,r_k;n-1}$ by
ignoring the last digit of $\pi$. For $\pi \in
I_{r_1,\dots,r_k;n}$ with $\pi(n)=k^{\vec{v}}$ and
$\pi(k)=n^{(r_1,\dots,r_k)-\vec{v}}$, define $\pi'' \in
I_{r_1,\dots,r_k;n-2}$ in the following way: Write $\pi$ in its
complete notation, i.e. as a matrix of two rows. The first row of
$\pi''$ is $(1,2,\dots,n-2)$ while the second row is obtained from
the second row of $\pi$ by ignoring the digits $n$ and $k$, and
the other digits are placed in an order preserving way with
respect to the second row of $\pi$. Here is an explicit formula
for the map $\pi \mapsto \pi''$.

$$\pi'' (i)=\left\{
\begin{array}{ccc}
\pi (i)      & & 1 \leq i < k\ {\quad {\rm and} \qquad} \pi(i) < k \\
\pi (i)-1    & & 1 \leq i < k\ {\quad {\rm and} \qquad} \pi(i) > k \\
\pi (i-1)    & &   k \leq i < n \ {\quad {\rm and} \qquad} \pi(i) < k \\
\pi (i-1)-1  & &   k \leq i < n \ {\quad {\rm and} \qquad} \pi(i) > k \\
\end{array}\right.$$

Note that the map $\pi \mapsto \pi'$ is a bijection from the set
$$\{\pi \in I_{r_1,\dots,r_k;n} \mid \pi(n)=n^{\vec{v}}\} \quad (\vec{v}\ {\rm fixed})$$
to $I_{r_1,\cdots,r_k;n-1}$, while $\pi \mapsto \pi''$ is a
bijection from the set
$$\{\pi \in I_{r_1,\dots,r_k;n} \mid
\pi(n)=k^{\vec{v}}\} \quad (\vec{v}\ {\rm fixed})$$ to
$I_{r_1,\dots,r_k;n-2}$.

For any $\vec{v}$, if $\pi(n)=n^{\vec{v}}$ then:
$$\fix(\pi)=\fix(\pi')+1,$$ $$\exc_A(\pi)=\exc_A(\pi').$$
Since $\vec{v}$ satisfies $2\vec{v}=0$, we have two cases. If
$\vec{v}=\vec{0}$, then:
$$\csum(\pi)=\csum(\pi'). $$
On the other hand, if $\vec{v}=(z_1,\dots,z_k)\neq \vec{0}$, then
there is some $m$, $1 \leq m \leq k$, such that $z_m\neq 0$, and
$z_i=0$ for all $1 \leq i <m$. In this case, $\vec{v}$ contributes
$z_m \cdot \frac{r}{r_m}$ to ${\rm csum}(\pi)$. But since $\pi$ is
an involution, we have $z_m =\frac{r_m}{2}$, and thus we have:
$$\csum(\pi)=\csum(\pi')+\frac{r}{2}. $$

If  $\pi(n)=t^{\vec{v_1}}$, then the parameters satisfy
$$\fix(\pi)=\fix(\pi''),$$
$$\exc_A(\pi)=\exc_A(\pi'')+\delta_{\vec{v},\vec{0}}.$$
where $\delta_{\vec{v}_1,\vec{v}_2}$ is a generalized Kronecker
Delta:
$$\delta_{\vec{v}_1,\vec{v}_2}=
\left\{ \begin{array}{cc} 1 & \vec{v}_1=\vec{v}_2 \\ 0 & \vec{v}_1
\neq \vec{v}_2
\end{array} \right..$$
Note that $\pi(t)=n^{\vec{v_2}}$, where $\vec{v}_1 +\vec{v}_2=0$.
Again, we have two cases. If $\vec{v}_1=\vec{0}$, then
$\vec{v}_2=\vec{0}$, and:
$$\csum(\pi)=\csum(\pi'). $$
On the other hand, if $\vec{v}_1=(z_1,\dots,z_k)\neq \vec{0}$, then
there is some $m$, $1 \leq m \leq k$, such that $z_m\neq 0$, and
$z_i=0$ for all $1 \leq i <m$. Since $\vec{v}_1 +\vec{v}_2=0$, we
have that $\vec{v}_2=(z'_1,\dots,z'_k)\neq \vec{0}$ with $z'_m\neq
0$, $z'_i=0$ for all $1 \leq i <m$, and $z_m+z'_m=r_m$. Now,
$\vec{v}_1$ contributes $z_m \cdot \frac{r}{r_m}$ to ${\rm
csum}(\pi)$ while $\vec{v}_2$ contributes $z'_m \cdot \frac{r}{r_m}$
to ${\rm csum}(\pi)$. Hence, their total contribution is:
$$z_m \cdot \frac{r}{r_m}+z'_m \cdot \frac{r}{r_m}=(z_m+z'_m)\cdot \frac{r}{r_m}=r_m \cdot \frac{r}{r_m}=r.$$
Thus we have:
$$\csum(\pi)=\csum(\pi'')+r(1-\delta_{\vec{v},\vec{0}}).$$

Now define $$\epsilon =\#\{r_i \ |\ 1\leq i \leq k, \quad r_i = 0
\pmod 2\}.$$ Define also
$$\mu=\mu_{r_1,\dots,r_k} = \left\{ \begin{array}{cc}  1 +2^\epsilon w^{\frac{r}{2}} & \epsilon\neq 0 \\ 1 & \epsilon=0 \end{array} \right..$$
The above consideration gives the following recurrence formula:
\begin{eqnarray*}
f_{r_1,\dots,r_k;n}(u,v,w) & = & u\mu f_{r_1,\dots,r_k;n-1}(u,v,w) \\
                    & & +(n-1)(v+(r-1)w^{r})f_{r_1,\dots,r_k;n-2}(u,v,w),\quad
                    n\geq1.
\end{eqnarray*}

Using the same technique used in \cite{BGM}, we get the following
explicit formula:
\begin{cor}
The polynomial $f_{r_1,\dots,r_k;n}(u,v,w)$ is given by

\begin{equation}\label{pol1}\sum_{j=n/2}^n
(n-j)!\binom{n}{n-j,n-j,2j-n}\frac{u^{2j-n}(v+(r-1)w^r)^{n-j}}{2^{n-j}}\mu^{2j-n}.\end{equation}
\end{cor}

If we substitute $w=1$ and compute the coefficient of
$u^mv^{\ell}$ in Formula (\ref{pol1}), we get the following
result:

\begin{cor}  The number of involutions in
$\mulgrn$ with exactly $m$ absolute fixed points and ${\rm exc}_A
(\pi)=\ell$ is given by
$$(\frac{n-m}{2})!(r-1)^{\frac{n-m}{2}-\ell}\binom{n}{\frac{n-m}{2},m,\frac{n-m}{2}-\ell,\ell}\frac{(1+2^{\epsilon})^{1-y}}{2^{\frac{n-m}{2}}},
$$ where $y \in \{0,1\}$ and $y \equiv r
\pmod 2$.
\end{cor}

We turn now to the computation of the number of involutions with a
fixed number of excedances. We do this by substituting $u=1$ and
$v=w^r$ in Formula (\ref{pol1}).

\begin{cor}\label{exccolor of grn}
The number of involutions  $\pi \in \mulgrn$ with $\exc(\pi)=m$ is:
\begin{equation*}\left\{\begin{array}{cc}
y!\binom{n}{y,\ y,\ n-2y}{(\frac{r}{2})}^y & \qquad r \equiv 1 \pmod 2 \\
\sumlim_{j=\frac{n}{2}}^{n}{(n-j)!\binom{n}{n-j,\ n-j,\ j-y,\
y-n+j}{(\frac{r}{2})}^{n-j}2^{\epsilon(y-n+j)}} & \qquad r \equiv 0
\pmod 2
\end{array} \right.,
\end{equation*}
where $y=\frac{m}{r}$.
\end{cor}

\section*{Acknowledgments}

The authors wish to thank Robert Schwartz for associating the
Dynkin-type diagram to the multi-colored permutation group. We
also thank Toufik Mansour for fruitful discussions.

\end{document}